\documentclass{amsart}

\usepackage{amsmath,amsthm}
\usepackage{amsfonts,amssymb}

\newtheorem{thm}{Theorem}[section]

\newtheorem{lem}[thm]{Lemma}

\theoremstyle{remark}

\def\a{{\alpha}} 
 \def\b{{\beta}}

 \def\l{{\lambda}}
 
 \def\om{{\omega}}
 \def\s{{\sigma}}
 \def\la{{\langle}}
 \def\ra{{\rangle}}

 \def\NN{{\mathbb N}}
 \def\CH{{\mathcal H}}

 \def\CP{{\mathcal P}}
 \def\CU{{\mathcal U}}
 \def\CV{{\mathcal V}}
 \def\RR{{\mathbb R}}

 \def\sspan{\operatorname{span}}

\begin{document}

\title[Orthogonal Polynomials and pde on unit ball]
{Orthogonal Polynomials and Partial Differential Equations on
the Unit Ball }

\author{Miguel Pi\~{n}ar}
\address{Department of Applied Mathematics\\
University of Granada\\
Granada 18071, Spain.}
\email{mpinar@ugr.es}
\author{Yuan Xu}
\address{Department of Mathematics\\ University of Oregon\\
    Eugene, Oregon 97403-1222.}\email{yuan@math.uoregon.edu}

\date{\today}
\keywords{pde, orthogonal polynomials, several variables, unit ball}
\subjclass{33C50, 33E30, 42C05}
\thanks{Partially supported by Ministerio de Ciencia y Tecnolog\'{\i}a
        (MCYT) of Spain and by the European Regional Development Fund
        (ERDF) through the grant MTM 2005--08648--C02--02, and Junta de
        Andaluc\'{\i}a, Grupo de Investigaci\'on FQM 0229.
        The work of the second author was supported in part by NSF 
        Grant DMS-0604056}

\begin{abstract}
Orthogonal polynomials of degree $n$ with respect to the weight function 
$W_\mu(x) = (1-\|x\|^2)^\mu$ on the unit ball in $\RR^d$ are known to 
satisfy the partial differential equation
$$
  \left[ \Delta  -  \la x, \nabla \ra^2 - 
         (2 \mu +d) \la x, \nabla \ra \right ] P =  -n(n+2 \mu+d) P
$$
for $\mu > -1$. The singular case of $\mu = -1,-2, \ldots$ is studied
in this paper. Explicit polynomial solutions are constructed
and the equation for $\nu = -2,-3,\ldots$ is shown to have complete 
polynomial solutions if the dimension $d$ is odd. The orthogonality of 
the solution is also discussed. 
\end{abstract}

\maketitle

\section{Introduction}\label{Introduction}
\setcounter{equation}{0}

Let $W_\mu(x) = (1-\|x\|^2)^\mu$. For $\mu > -1$, $W_\mu$ 
is integrable on the unit ball $B^d := \{x: \|x\| \le 1\}$ and 
orthogonal polynomials with respect to $W_\mu$ are well
defined. We denote by $\CV_n^d(W_\mu)$ the space of 
orthogonal polynomials of degree $n$ with respect to $W_\mu$ 
on $B^d$. It is well known (cf. \cite{DX}) that $P \in \CV_n^d(W_\mu)$ 
satisfies the partial differential equation
\begin{align}  \label{pde}
    L_\mu P:= &  \left[ \Delta  -  \la x, \nabla \ra^2 - 
         (2 \mu +d) \la x, \nabla \ra \right ] P = \l_n^{(\mu)} P, \quad 
      \hbox{with} \\
     \l_n^{(\mu)}:= &  -n(n+2 \mu+d), \notag
\end{align}
where $\Delta$ denotes the usual Laplacian, $\nabla = (\partial_1,
\ldots, \partial_d)$ with $\partial_i =  \frac{\partial}{\partial x_i}$ so
that $\la x, \nabla \ra  = \sum_{i=1}^d x_i \frac{\partial}{\partial x_i}$. 
In other words, the 
elements of $\CV_n^d(W_\mu)$ are eigenfunctions of $L_\mu$ 
with eigenvalues $\l_n^{(\mu)}$. Since $\cup_{n \ge 0} \CV_n^d(W_\mu)$ 
contains a basis for all polynomials, we say that $L_\mu$ has  
complete polynomial solutions for $\mu > -1$. 

The purpose of this paper is to study  equation \eqref{pde} for 
$\mu = -1, -2, \dots$. Our study starts from the realization that
orthogonal polynomials with respect to the inner product 
$$
   \la f, g \ra_{-1} =  \l \int_{B^d} \nabla f (x) \cdot \nabla g(x) dx + 
    \int_{S^{d-1}} f(x)g(x) d\omega, \qquad \l > 0, 
$$
which were studied recently in \cite{X07}, satisfy \eqref{pde} for 
$\mu = -1$ (see Section 2).  
To verify it, explicit orthogonal structure from \cite{X07} is used. 
It turns out that it is possible to construct explicit polynomial solutions 
for the equation \eqref{pde} when $- \mu \in \NN$. However, 
somewhat surprisingly, \eqref{pde} has a complete set of polynomial 
solutions only when the dimension $d$ is odd, and our construction 
shows why this is so. A natural question is whether the polynomial 
solutions of \eqref{pde} remain orthogonal for $\mu = -2, -3, \ldots$ 
when $d$ is odd. We do not know the answer to this question. We are, 
however, able to show that a family of polynomials closely related to 
the solutions of the equation \eqref{pde} for $k = -2$ are orthogonal
with respect to 
$$
  \la f, g \ra_{-2}
   =  \l \int_{B^d} \Delta f (x) \Delta g(x) dx + 
    \int_{S^{d-1}} f(x)g(x) d\omega, \qquad \l > 0.   
$$
Furthermore, for each $n \in \NN$, the polynomial solutions of 
\eqref{pde} for $k = -2$ are orthogonal with respect to an inner product 
that is a modification of $\la \cdot, \cdot \ra_{-2}$ but depends on $n$.

Historically, only the case $d =1$ of such a problem has been thoroughly
studied. In the case of $d =1$, the equation \eqref{pde} becomes the 
ordinary differential equation satisfied by the Gegenbauer polynomials. In
this case, the problem of negative indices has been studied by several 
authors, we refer to \cite{AAR, APPM, APP, JKL} and the references
therein. For $d =2$, equation \eqref{pde} is classical and can be 
traced back to Hermite. It is one of five second order partial differential 
equations with polynomial coefficients that have orthogonal polynomials
solutions with respect to a positive definite moment functional, as 
classified by Krall and Scheffer in \cite{KS}. We note that our result for 
$\mu =-1$ does not conflict with the result of \cite{KS}, since 
$\la f, g \ra_{-1}$ does not come from a moment functional. For $d =2$, 
the case $- \mu \notin \NN$ is considered in \cite{LL2}, where
the fact that \eqref{pde} does not have complete polynomial solutions 
for $- \mu \in \NN$ is also observed.  Finally, we should mention \cite{LL}
and \cite{X06} for relevant results. 

The paper is organized as follows. We collect background materials
in the following section. The solutions of \eqref{pde} for $\- \mu \in \NN$ 
are constructed in Section 3. Orthogonal polynomials with respect to 
$\la\cdot, \cdot \ra_{-2}$ are discussed in Section 4.

\section{Background and Preliminary}\label{BP}
\setcounter{equation}{0}

Let $d\omega$ denote the Lebesgue measure on $S^{d-1}:=\{x:  \|x\| =1\}$
and denote the area of $S^{d-1}$ by $\omega_{d}$,  
$
    \om_{d} := \int_{S^{d-1}} d\om = 2 \pi^{d/2}/\Gamma(d/2).
$
We start with the spherical harmonics, which are 
homogeneous polynomials satisfying the equation $\Delta P =0$. Let 
$\CH_n^d$ denote the space of homogeneous harmonic polynomials of 
degree $n$. It is well known that 
$$
\dim \CH_n^d = \binom{n+d-1}{d-1}-\binom{n+d-3}{d-1} := \s_n.
$$ 
The restriction of $Y \in \CH_n^d$ on $S^{d-1}$ are the spherical harmonics,
which are orthogonal on $S^{d-1}$. Throughout this paper, we use the 
notation $\{Y_\nu^n: 1 \le \nu \le \sigma_n\}$ to denote an orthonormal basis 
for  $\CH_n^d$, that is,  
\begin{equation} \label{eq:harmonics}
\frac{1}{\omega_{d}} \int_{S^{d-1}} Y_\mu^n(x') Y_\nu^m(x') d\omega(x')
    = \delta_{\mu,\nu} \delta_{n,m} , \qquad 1 \le \mu,\nu\le \sigma_n.
\end{equation}
We will also denote by $\Pi_n^d$ the space of polynomials of degree $n$
in $d$ variables and by $\CP_n^d$ the space of homogeneous polynomials
of degree $n$. It is well known that $\dim \CP_n^d = \binom{n+d-1}{d-1}$. 

For orthogonal polynomials on the unit ball, we start with the inner product
$$
  \la f, g \ra_\mu =  c_\mu \int_{B^d} f (x) g(x)   
       (1-\|x\|^2)^\mu dx, \qquad \mu > -1. 
$$
where $c_\mu$ is the normalization constant of $W_\mu$. Let 
$\CV_n^d(W_\mu)$ denote the space of orthogonal polynomials of 
degree $n$. A mutually orthogonal basis for $\CV_n^d(W_\mu)$  is 
given by (\cite{DX}) 
\begin{equation}\label{eq:Wmu-basis}
 P_{j,\nu}^n(W_\mu; x) = P_j^{(\mu, n-2j+\frac{d-2}{2})}(2\|x\|^2 -1) 
   Y_{\nu}^{n-2j}(x), \quad 0 \le j \le n/2,
\end{equation}
where $P_j^{(\a,\b)}$ denotes the Jacobi polynomial of degree $j$, 
which is orthogonal with respect to $(1-x)^\alpha(1+x)^\beta$ on $[-1,1]$, 
and $\{Y_\nu^{n-2j}: 1 \le \nu \le \sigma_{n-2j}\}$ is an orthonormal basis 
for $\CH_{n-2j}^d$. 

We will need orthogonal polynomials for two other inner products on
the ball. One is $\la\cdot, \cdot \ra_{-1}$ in the introduction, which we
normalize as  
$$
  \la f, g\ra_{-1} := \frac{\l }{\om_d} 
      \int_{B^d} \nabla f(x) \cdot \nabla g(x)dx 
     + \frac{1}{\om_d} \int_{S^{d-1}} f(x) g(x) d\om_d, \quad \l > 0,
$$
so that  $\la 1, 1 \ra_{-1} =1$.  Let $\CV_n^d(W_{-1})$ denote the 
space of orthogonal polynomials with respect to $\la 1, 1\ra_{-1}$. 
Then, as shown in \cite{X07},  
\begin{equation}\label{V{-1}}
   \CV_n^d(W_{-1}) = \CH_n^d \oplus (1- \|x\|^2) \CV_{n-2}^d(W_1). 
\end{equation}
The second inner product is motivated by an application in 
numerical solution of Poisson equations (see \cite{A}) and it is 
defined by
\begin{equation}\label{eq:InnerDelta}
  \la f, g \ra_\Delta =  \a_d \int_{B^d} \Delta[(1-\|x\|^2) f (x)]  
           \Delta[(1-\|x\|^2) g (x)] dx, 
\end{equation}
where $a_d = 1/(4 d^2 {\rm vol} (B^d))$ so that $\la 1,1\ra_\Delta =1$.
Let $\CV_n^d(\Delta)$ denote the space of orthogonal polynomials 
with respect to $\la \cdot, \cdot\ra_{\Delta}$. Then, as shown in \cite{X06}, 
\begin{equation}\label{VDelta}
\CV_n^d(\Delta)  = \CH_n^d \oplus (1- \|x\|^2) \CV_{n-2}^d(W_2), 
\end{equation}
We refer to \cite{X06,X07} for further properties of orthogonal 
polynomials with respect to these two inner products. Note that an
explicit orthonormal basis for either orthogonal polynomial subspace
can be easily deduced from \eqref{eq:Wmu-basis}.

In the above and throughout this paper, we write $g(x) V_n = 
\{g(x) p(x): p \in V_n\}$, where $V_n$ is either  the subspace 
$\CH_n^d$ of harmonic polynomials or the subspace 
$\CV_n^d(W_\mu)$ of orthogonal polynomials. Furthermore, if
$n < 0$, then we define $\CH_n^d = \emptyset$ and $\CV_n^d = \emptyset$.

\section{Polynomial Solutions of the Partial Differential 
Equation}\label{OP-PDE}
\setcounter{equation}{0}

We start with two lemmas on the operator $L_\mu f$. 

\begin{lem}  \label{Lem1}
For $ \mu \in \RR$, 
\begin{equation} \label{recursive}
   L_\mu  [(1-\|x\|^2)P]   = 4 ( \mu +1) P   + 
    (1-\|x\|^2)  \left[ -  (4 (\mu+1) +2 d ) P  + L_{\mu +2} P \right] .  
\end{equation} 
\end{lem} 

\begin{proof} 
The proof follows from a straightforward computation. The 
essential ingredients are 
\begin{align*}
   \Delta[ (1-\|x\|^2)P] = & -2 d P - 4 \la x, \nabla \ra P+ (1-\|x|^2)P, \\
   \la x, \nabla \ra [ (1-\|x\|^2) P]   = & -2 P  + 2 (1-\|x\|^2) P  +  
       (1-\|x\|^2) \la x, \nabla \ra P + (1-\|x|^2)P, \\ 
 \la x, \nabla \ra^2 [(1-\|x\|^2) P] = & -4 P  - 4 \la x, \nabla \ra P  
        + 4 (1-\|x\|^2) P   \\ &
           +  4 (1-\|x\|^2) \la x, \nabla \ra P  + (1-\|x|^2)\la x, \nabla \ra^2 P, 
\end{align*}
as can be easily verified. 
\end{proof} 

\begin{lem} \label{Lem2} 
For all $\mu \in \RR$ and $n \ge 0$, $L_\mu Y =   \l_n^{(\mu)} Y$ for 
$Y \in \CH_n^d$. 
\end{lem}

\begin{proof} 
If $Y \in \CH_n^d$, then $\Delta Y =0$. Since $Y$ is homogeneous, 
Euler's formula shows that $\la x, \nabla \ra Y = n Y$. The stated
result follows immediately from these two facts. 
\end{proof} 

Our first result shows that orthogonal polynomials with respect to
$\la\cdot,\cdot \ra_{-1}$ satisfy the equation \eqref{pde} with
$\mu = -1$. 

\begin{thm} \label{thm1}
Elements of $\CV_n^d(W_{-1})$ satisfy $L_{-1} P =  \l_n^{(-1)} P$. 
In particular, the eigenfunctions of the operator $L_{-1}$ consist of 
a complete polynomial  basis. 
\end{thm} 

\begin{proof}
Setting $\mu = -1$, equation \eqref{recursive} shows that 
$$
 L_{-1}   \left[ (1-\|x\|^2) P\right] = (1-\|x\|^2)  \left[ - 2 d P +
    L_1 P \right].
$$
Recall that elements of $\CV_{n-2}^d(W_1)$  satisfy the equation
\eqref{pde} with $\mu =1$. It follows that for $P \in \CV_{n-2}^d(W_1)$, 
\begin{align*}
 L_{-1} \left[ (1-\|x\|^2) P\right] = &  (1-\|x\|^2) \left[ - 2 d P -  
       (n-2) (n + 2\mu + d -2) P \right] \\
        = & -n (n+d -2)  (1-\|x\|^2) P. 
\end{align*}
By \eqref{V{-1}}, this and Lemma \ref{Lem2} show that all elements
in $\CV_n^d(W_{-1})$ satisfy \eqref{pde} for $\mu = -1$. 
\end{proof}

For $d =2$, a complete characterization of second order linear 
partial differential equations that have polynomial coefficients and 
orthogonal polynomials as eigenfunctions is given in \cite{KS},
where the orthogonality is given in terms of a moment functional.
If the moment functional is positive definite, then there are only
five families of equations, the equation \eqref{pde} is one of them
when $\mu > -1$. It should be remarked that Theorem \ref{thm1}
does not contradict the characterization of \cite{KS}, since the
inner product $\la\cdot, \cdot \ra_{-1}$ cannot be expressed
as a moment functional. Indeed, it is easy to see that 
$\la x f, g \ra_{-1} \ne \la  f, x g \ra_{-1}$ in general.

The above theorem prompts us to consider solutions of \eqref{pde} 
for $\mu = -2,-3,\ldots$ and leads to our second theorem. 
For $\mu = -k$, $k = 2, 3, \ldots$, we define 
$$%\begin{equation}
   \CU_n^d(W_{-k}) := \CH^d_n \cup \left(\bigcup_{j=1}^{k-1}
       \left[  \sum_{\nu=0}^{j} a_{j,\nu}^n
      (1-\|x\|^2)^\nu \right] \CH_{n-2j}^d \right) \cup
        (1-\|x\|^2)^k \CV_{n-2k}^d(W_k),
$$%\end{equation}
where, for  $1 \le j \le k-1$,
\begin{equation} \label{a_jnu}
  a_{j,\nu}^n := 
    \frac{(-1)^{j - \nu} j! (-k+1)_j (n-j-k+\frac{d}{2})_\nu}
    {\nu! ( j - \nu )!  (-k+1)_\nu (n-j-k+\frac{d}{2})_j}, \quad 0 \le \nu\le j, 
\end{equation}
in which $(a)_m : = a(a+1) \ldots (a+m-1)$ is the shifted factorial. 
We note that $a_{j,j}^n = 1$. It is easy to see that $a_{j,\nu}^n$ are well 
defined if $n - j - k +\nu + \frac{d}{2}$ is not zero. 

\begin{thm} \label{thm2}
If $\mu = - k$ and $k = 2, 3, \ldots$,  then the polynomials in
$\CU_n^d(W_{-k})$ satisfy equation \eqref{pde}; that is,
$L_{-k} P = \lambda_n^{(-k)} P$ for $P \in \CU_n^d(W_{-k})$. 
Furthermore, 
$$
    \dim \CU_n^d = \dim \CP_n^d, \qquad 
       \hbox{if \quad $n - 2 k -1 +\frac{d}{2} \ne 0$. }
$$
In particular,  the operator $L_{-k}$  has a complete polynomial basis 
of eigenfunctions if the dimension $d$ is odd. 
\end{thm}

\begin{proof}
We observe that the relation \eqref{recursive} can be used recursively,
which yields 
\begin{align} \label{recursive2}
    L_{\mu}(1-\|x\|^2)^k =\, & 4 k (\mu+k) (1-\|x\|^2)^{k-1}  \\
         & + (1-\|x\|^2)^k \left[ - (4k(\mu+k) + 2 k d) + L_{\mu+2k}\right], \quad
            k \in \NN, \notag
\end{align}
as can be easily verified by induction. Hence, for $Q \in 
\CV_{n-2k}^d(W_k)$, which satisfies $L_k Q = \lambda_{n-2k}^{(k)}Q$,  
we have 
\begin{align*}
 L_{-k} [(1-\|x\|^2)^k Q(x)] &=  (1-\|x\|^2)^k \left[ - 2 k d Q(x) + L_k Q(x) \right] \\
     & =  (1-\|x\|^2)^k \left[ - 2 k d  - (n-2 k) (n-2 k + 2k +d) \right] Q(x) \\
     & = -n (n-2k+d)  (1-\|x\|^2)^k Q(x),
\end{align*}
so that elements in  $(1-\|x\|^2)^k \CV_{n-2k}^d(W_k)$ satisfy 
$L_{-k} P = \lambda_n^{(-k)}P$. Furthermore, by Lemma \ref{Lem2},
elements of $\CH_n^d$ are also eigenfunctions of $L_{-k}$ with
eigenvalue $\l_n^{(-k)}$. Let now
$$
   Z_j(x) =  (1-\|x\|^2)^j Y(x) +\sum_{\nu=0}^{j-1} a_{j,\nu}^n
       (1-\|x\|^2)^\nu Y(x),    \qquad Y \in \CH_{n-2j}^d. 
$$
We want to show that $Z_j$ satisfies the same equation. Applying 
\eqref{recursive2} shows that 
\begin{align*}
 L_{-k} Z_j(x) = & - n(n-2 k + d) (1-\|x\|^2)^j Y(x)+ \sum_{\nu=0}^j 
     \left[4 (\nu +1)(-k + \nu+1)a_{j,\nu+1}^{n} \right. \\
   & \left. - (n-2j+ 2 \nu)(n-2j + 2 \nu - 2 k +2) a_{j,\nu}^{n} \right]
       (1-\|x\|^2)^\nu Y(x),
\end{align*}  
so that the equation $L_{-k} Z_j = \lambda_n^{(-k)} Z_j$ implies, 
upon comparing the coefficients, that $a_{j,\nu}^{(n)}$ satisfy
the recursive relation
$$
    a_{j,\nu}^{(n)}= a_{j,\nu+1}^{(n)} \frac{-2 (\nu+1)(-k+\nu+1)}{(j-\nu) 
       (2n- 2j + 2\nu-2k+d)}. 
$$
Setting $a_{j,j}^{(n)} =1$, we can then determine $a_{j,\nu}^{(n)}$
recursively. The result is \eqref{a_jnu}. This shows that all elements
in $\CU_n^d(W_{-k})$ satisfy the equation $L_{-k} P = \l_n^{(-k)} P$. 

If $a_{j,\nu}^n$ are all finite, then it follows from 
$\dim \CP_n^d = \dim \CH_{n}^d + \dim \CP_{n-2}^d$ that 
$$
\dim \CU_n^d(W_{-k}) = \dim \CH_n^d + \ldots + \dim \CH_{n-2k+2}^d 
         + \dim \CP_{n-2k}^d = \dim \CP_n^d.
$$ 
The explicit formula \eqref{a_jnu} shows that $a_{j,\nu}^n$ is not finite
only if $n - j - k +\nu + d/2$ is zero, which can happen only if 
$d$ is even and $n < j+k - \nu- d/2$. Since $0 \le \nu \le j-1$ and 
$j \le k-1$, this could happen only if $n - 2k-1+d/2$ is zero or a 
negative integer.
\end{proof} 
 
Some remarks are in order. If $d$ is even then, as shown in the proof, 
$\CU_n^d(W_{-k})$ is not well defined if $n - j - k +d/2 =0$. The first 
time this happens is when $d =2$ and $k =2$. As $ 1 \le j \le k-1$, this
happens when $j=1$ and $n = 2$.  In fact, for $k =2$, we have 
\begin{equation*} %\label{CU2}
 \CU_n^2(W_{-2}) = \CH_n^2 \cup \left[ (1-\|x\|^2) + \frac{1}{n-2}\right]
    \CH_{n-2}^2 \cup (1-\|x|^2)^2 \CV_{n-4}^2(W_2),
\end{equation*}
so that $\CU_2^d(W_{-2})$ is not well defined. Furthermore, since 
$\l_2^{(-2)} =0$, $\CH_2 = \sspan\{x_1^2 - x_2^2, 2 x_1 x_2\}$ 
and $L_\mu 1 =0$, we see that the linearly independent eigenfunctions
with respect to the eigenvalue $\l_2^{(-2)} =0$ are  $1, x_1^2-x_2^2,
2 x_1x_2$. However, $1$ is also the eigenfunction of  $\l_0^{(-k)}$, 
we see that $L_{-2}$ does not have a complete polynomial  basis as
eigenfunctions and it misses exactly one polynomial of degree $2$. 

The same analysis also works for  $k \ge 2$. For example, for $k =4$
and $d =2$,  
\begin{align*}
 \CU_n^2(W_{-4}) =   \CH_n^2 & \cup \left[ (1-\|x\|^2) + \frac{3}{n-4}\right]
     \CH_{n-2}^2    \\
     &\cup \left[(1-\|x|^2)^2 + \frac{4}{n-4} (1-\|x\|^2) 
         + \frac{8}{(n-4)(n-5)} \right] \CH_{n-4}^2   \\
      & \cup \left[ (1-\|x|^2)^3 + \frac{3}{n-4} (1-\|x\|^2)^2 
         + \frac{6}{(n-4)(n-5)}  (1-\|x\|^2) \right. \\
       & \left.   \qquad + \frac{6}{(n-4)(n-5)(n-6)} \right] \CH_{n-6}^2  
           \cup (1-\|x|^2)^2 \CV_{n-8}^d(W_2), 
\end{align*}
which shows clearly that the operator $L_{-4}$ has a complete basis 
of polynomials of degree exactly $n$ as eigenfunctions with respect to 
the eigenvalue $\l_n^{(-k)}$ only if  $n > 6$. 

In general, for even dimension $d$, the operator $L_-k$ has a complete 
basis of polynomials of degree exactly $n$ as eigenfunctions only if 
$n > 2k-1+d/2$.  As a consequence, the eigenfunctions of the operator
$L_{-k}$ do not contain a complete basis of polynomials if $d$ is even.

\section{Orthogonal polynomials associated with Laplacian}
\setcounter{equation}{0}

For all $\mu \in \RR$, the function $W_\mu$ is a symmetric factor for
the differential operator $L_\mu$ as it is easy to verify that 
$L_\mu W_\mu = (W_\mu L_\mu)^*$ with
$$
(W_\mu L_\mu)^* := \Delta W_\mu + (d+2 \mu) \la \nabla, x\ra 
  W_\mu -  \la \nabla, x\ra^2  W_\mu. 
$$
The verification is formal and holds for all $\mu \in \RR$. The 
polynomial eigenfunctions of $L_\mu$ are orthogonal with 
respect to $W_\mu$ for $\mu > -1$.  According to Theorem \ref{thm1}, 
the polynomial eigenfunctions of $L_{-1}$ form a complete orthogonal 
basis with respect to $\la\cdot, \cdot \ra_{-1}$.  It is natural to ask if 
the same holds true for $L_{-k}$ with $k = 2, 3, \ldots$  for $d$ being 
odd. Below we consider the case $k =2$. 

Recall that  in the case of $k = 2$, the eigenfunctions of  
$L_{-2}$ are given by 
\begin{equation}\label{CU2}
 \CU_n^d(W_{-2}) = \CH_n^d \cup \left[ (1-\|x\|^2) + 
 \frac{1}{n-3+\frac{d}{2}}\right]
  \CH_{n-2}^d \cup (1-\|x\|^2)^2 \CV_{n-4}^d(W_2). 
\end{equation}
We first consider the following closely related class of functions 
\begin{equation} \label{CV2}
 \CV_n^d(W_{-2}): = \CH_n^d \cup  (1-\|x\|^2)\CH_{n-2}^d \cup
       (1-\|x\|^2)^2 \CV_{n-4}^d(W_2). 
\end{equation}
The difference between $\CV_n^d(W_{-2})$ and $\CV_n^d(W_{-2})$ 
are only in the middle term. Notice that $\CV_n^d(W_{-2})$ are defined
for all $n$ and $d$. It turns out that it is the space of orthogonal 
polynomials with respect to the inner product
$$
   \la f, g \ra_{-2}: = \frac{\l}{\om_d} \int_{B^d} \Delta f(x) \Delta g(x) dx + 
   \frac{1}{\om_d}   \int_{S^{d-1}} f(x) g(x) d\om, \quad \l > 0,
$$
which is normalized so that $\la 1,1\ra_{-2} =1$. It is easy to see that 
this inner product is well-defined and it is positive definite.

\begin{thm} \label{thm3}
The elements of $\CV_n^d(W_{-2})$ are orthogonal polynomials 
with respect to $\la \cdot, \cdot \ra_{-2}$. Moreover, they contain an
orthonormal basis; in other words,
$$
    \CV_n^d(W_{-2})  = \CH_n^d  \oplus  (1-\|x\|^2)\CH_{n-2}^d \oplus
       (1-\|x\|^2)^2 \CV_{n-4}^d(W_2).
$$
\end{thm}

\begin{proof}
First we make an observation that, by \eqref{VDelta}, 
$$
   \CV_n^d (W_{-2}) =  \CH_n^d \cup (1-\|x\|^2) \CV_{n-2}^d(\Delta). 
$$
We then prove that each of the two components in the right hand side
of the above expression are orthogonal to lower degree polynomials 
and they are mutually orthogonal with respect to $\la\cdot,\cdot\ra_{-2}$.
We consider the two cases separately. 

\medskip\noindent{\it Case 1.} Let $P=Y_n \in \CH_n^d$. Since $\Delta P =0$,
$\la P, g \ra_{-2}  =  \frac{1}{\om_d}   \int_{S^{d-1}} f(x) g(x) d\om$. If $g 
\in \Pi_{n-1}^d$, then $\la P, g \ra_{-2} =0$ by the orthogonality of $Y_n$. 
If $g\in (1-\|x\|^2)\CV_{n-2}^d(\Delta)$, then the restriction of $g$ on $S^{d-1}$ is 
zero, so that $\la P, g \ra_{-2} =0$. 

\medskip\noindent{\it Case 2.} Let $P(x) = (1-\|x\|^2) Q_{n-2}(x)$, where 
$Q_{n-2} \in \CV_{n-2}^d(\Delta)$. Then the restriction of $P$ on $S^{d-1}$ 
is zero, so that 
$$
\la P, g \ra_{-2}  =  \frac{\l}{\om} \int_{B^d} \Delta [ (1-\|x\|^2) Q_{n-2}(x)] 
    \Delta g(x) dx.
$$
Let $g \in \CV_m^d(W_{-2})$ and $ m < n$.  If $g \in \CH_m$, then $\Delta g
=0$ so that $\la P, g \ra_{-2} =0$. If $g \in (1-\|x\|^2) \CV_{m-2}^d(\Delta)$, 
then $g(x) = (1-\|x\|^2) h(x)$ for $h \in \CV_{m-2}^d(\Delta)$. Recall the 
definition of $\la f , g\ra_\Delta$ in \eqref{eq:InnerDelta}, we then have
$\la P, g \ra_{-2} = \la Q_{n-2}, h \ra_{\Delta} =0$ by the orthogonality of
$Q_{n-2} \in  \CV_{n-2}^d(\Delta)$.
\end{proof}

Recall that $\CV_n^d(W_{-2})$ and $\CU_n^d(W_{-2})$ differ only at the middle
component. We do not know if there is an inner product with respect to which 
$\CU_n^d(W_{-2})$ is orthogonal for {\it all} $n$.  For each fixed $n \in \NN$, 
however,  elements of $\CU_n^d(W_{-2})$ are orthogonal to lower degree 
polynomials with respect to the following inner product: 
\begin{align*}
\la f, g \ra_n^*: = & \frac{\l}{\om_d} \int_{B^d} \Delta f(x) \Delta g(x) dx + 
    \frac{1}{\om_d}   \int_{S^{d-1}} f(x) g(x) d\om \\
 & +  \frac{\mu}{\om_d} \int_{S^{d-1}} \frac{d}{dr} \left[r^{n-4-d} f(x)\right]
      \frac{d}{dr} \left[r^{n-4-d} g(x)\right] d\om,
\end{align*}
where $\frac{d}{dr}$ is the normal derivative, $\l>0$ and $\mu > 0$. It is 
easy to see that $\la f, g \ra_n^*$ is positive definite for $f, g \in \Pi_n^d$.

\begin{thm}
The elements of $\CV_n^d(W_{-2})$ are orthogonal polynomials 
with respect to the inner product $\la f, g \ra_n^*$.
\end{thm}

\begin{proof}
Notice that the first two terms in the $\la f, g \ra_n^*$ are exactly 
$\la f, g \ra_{-2}$.  For $x \in \RR^d$, write $x = r x '$, $r > 0$ and 
$x' \in S^{d-1}$. We assume that $g \in \Pi_{n-1}^d$ below. 
If $ P = Y_n \in \CH_n^d$, then $\frac{d}{dr} \left[r^{n-4-d} P(x)\right] 
= (2n-4-d) Y_n(x')$ on $S^{d-1}$. Hence $\la P, g \ra_n^* = 0$ as 
$\Delta P =0$ and $P \in CH_n^d$. If $P \in (1-\|x\|^2)^2 \CV_{n-4}^d(W_2)$, 
then $\frac{d}{dr} \left[r^{n-4-d} P(x)\right] =0$ on $S^{d-1}$ because of 
the $(1-r^2)^2$ term, so that $\la P, g \ra_n^* = \la P, g \ra_{-2} = 0$. 
Thus, by \eqref{CU2},  we are left with the case 
$P(x)= \left[(1-\|x\|^2) + \frac{1}{n-3+d/2}\right] Y_{n-2}(x)$, where
$Y_{n-2}  \in  \CH_{n-2}^d$. Since $Y_{n-2}$ is homogeneous,
$Y_{n-2}(x) = r^{n-2} Y_{n-2}(x')$, a simple computation shows that
$\frac{d}{dr} \left[r^{n-4-d} P(x)\right] =0$. Consequently, we again have
$\la P, g \ra_n^* = \la P, g \ra_{-2} = 0$. 
\end{proof}

These results indicate that there is a difference between the case $k = -1$
and $k =-2,-3, \ldots$. Not only the differential operator $L_{-k}$ has 
complete polynomial eigenfunctions only when $d$ is odd in the latter 
case, the eigenfunctions also may no longer be orthogonal with respect
to an inner product for all $n$. 

Another interesting question is if the main terms of $\CU_n^d(W_{-k})$ 
are orthogonal for $k \ge 3$. In other words, if
the spaces 
$$
  \CV_n^d(W_{-k}): = \CH^d_n \cup \left(\bigcup_{j=1}^{k-1}
       (1-\|x\|^2)^j \CH_{n-2j}^d \right) \cup
        (1-\|x\|^2)^k \CV_{n-2k}^d(W_k),
$$
are orthogonal with respect to a positive definite inner product, which 
will likely involve higher order derivatives, as indicated by our result for 
$k =1$ and $k=2$. At this point we do not see a way to find such an
inner product.

\medskip\noindent
{\it Acknowledgment.} The authors thank Jeong Keun Lee for 
communicating and discussing the results in \cite{LL2}.

\end{document}